\documentclass[11pt]{amsart}
\title{Applications of the Defect of a Finitely Presented Functor}
\author{Jeremy Russell}
\address{The College of New Jersey, Ewing, New Jersey}
\email{russelj1@tcnj.edu}
\usepackage{amsmath,amssymb,amsthm, amsfonts, xypic, graphicx}

\usepackage[hidelinks=true]{hyperref}
\usepackage[all]{xy}

\theoremstyle{definition}
\newtheorem{thm}{Theorem}
\newtheorem{prop}[thm]{Proposition}
\newtheorem{lem}[thm]{Lemma}
\newtheorem{cor}[thm]{Corollary}

\newenvironment{pf}{\paragraph{{\sc \dbf{Proof}}}}{\par\vspace{1cm}}
\theoremstyle{definition}
\newtheorem{defn}{\underline{Definition}}

\newtheorem*{thm*}{Theorem}

\newtheorem*{com}{Comment}
\newtheorem*{cor*}{Corollary}
\newtheorem*{lem*}{Lemma}
\newtheorem{ex}{Example}

\renewcommand{\qed}{\blacksquare}

\newcommand{\blank}{\hspace{0.05cm}\underline{\ \ }\hspace{0.1cm} }

\newcommand{\ZZ}{\mathbb{Z}}

\renewcommand{\mod}{\textsf{mod}}
\newcommand{\Mod}{\textsf{Mod}}
\newcommand{\A}{\mathcal{A}}
\newcommand{\B}{\mathcal{B}}

\renewcommand{\S}{\mathcal{S}}
\newcommand{\Y}{\textsf{Y}}

\newcommand{\sub}{\subseteq}

\newcommand{\fp}{\textsf{fp}}

\newcommand{\dbf}[1]{\textbf{\text{#1}}}

\newcommand{\dia}[1]{\[\xymatrix{#1 }\]}
\usepackage{pgfpages}
\usepackage{tikz}
\usetikzlibrary{matrix,arrows,snakes,backgrounds}

\renewcommand{\:}{\colon}

\newcommand{\eps}{\varepsilon}
\newcommand{\To}{\longrightarrow}

\newcommand{\lex}{\textsf{Lex}}

\newcommand{\ab}{\textsf{Ab}}

\newcommand{\coker}{\textsf{Coker}}
\renewcommand{\ker}{\textsf{Ker}}

\newcommand{\ext}{\textsf{Ext}}

\newcommand{\tr}{\textrm{Tr}}

\renewcommand{\u}{\Omega}

\newcommand{\nat}{\textsf{Nat}}
\renewcommand{\hom}{\textsf{Hom}}

\newcommand{\testarr}[4]{(m-#1-#2)edge node[auto]{}(m-#3-#4)}

\begin{document}
\maketitle

\begin{abstract}For an abelian category $\A$, the defect sequence 
$$0\To F_0\To F\overset{\varphi}{\To} \big(w(F),\blank\big)\To F_1\To 0$$ of a finitely presented functor is used to establish the CoYoneda Lemma.  An application of this result is the \fp-dual formula which states that for any covariant finitely presented functor $F$,  $F^*\cong \big(\blank, w(F)\big)$.  The defect sequence is shown to be isomorphic to both the double dual sequence $$0\To \ext^1(\tr F,\hom)\To F\To F^{**}\To \ext^2(\tr F,\hom)\To 0$$ and the injective stabilization sequence $$0\To \overline{F}\To F\To R^0F\To \tilde F\To 0$$ establishing the \fp-injective stabilization formula $\overline{F}\cong \ext^1(\tr F,\hom)$ for any finitely presented functor $F$.   The injectives of $\fp(\Mod(R),\ab)$ are used to  compute the left derived functors $L^k(\blank)^*$.  These functors are shown to detect certain short exact sequences.  \end{abstract}

\setcounter{tocdepth}{1}

\tableofcontents

\section{Introduction}

For any abelian category $\A$, the category of finitely presented functors $\fp(\A,\ab)$ consists of all functors $F\:\A\to \ab$ for which there exists an exact sequence $$(Y,\blank)\To (X,\blank)\To F\To 0$$  One of Auslander's major contributions to representation theory was the demonstration that one may study the category $\A$ by studying the category $\fp(\A,\ab)$ of finitely presented functors.  This originates in \cite{coh} and is continued in many subsequent works such as \cite{repdim}, \cite{fart}, and \cite{stab}.  The functorial techniques Auslander used to study finitely presented functors have become widespread in representation theory of algebras but have also been applied to different fields such as algebraic geometry and model theory.   For more information, the reader is referred to \cite{cohhart} and  \cite{psl}.

This paper focuses on applications of the defect which is an exact contravariant functor $$w\:\fp(\A,\ab)\To \A$$ constructed by Auslander in \cite{coh}.  Given a finitely presented functor $F$ and a presentation $$(Y,\blank)\To (X,\blank)\To F\To 0$$ the functor $w$ is completely determined by the exact sequence $$0\To w(F)\To X\To Y$$  The functor $w$ reveals information at a local level about the functor $F$ and from a more macroscopic point of view it reveals information about the category $\A$.   Auslander also associates to each finitely presented functor $F$, an exact sequence $$0\To F_0\To F\To \big(w(F),\blank\big)\To F_1\To 0$$  called the defect sequence.   Years later in \cite{stab}, Auslander and Bridger define for any additive functor $F\:\A\To \B$ an injective stabilization sequence $$0\To \overline{F}\To F\To R^0F\To \tilde F\To 0$$ This requires that $\A$ has enough injectives and that $\B$ is abelian.   Understanding the relationship between these two sequences is the central motivation behind studying the defect in more detail.  Our approach leads to the following main results.

\begin{lem*}[The CoYoneda Lemma]For any finitely presented functor $F$, $$\nat\big(F,\big(X,\blank)\big)\cong \big(X,w(F)\big)$$\end{lem*}

\begin{cor*}[\fp-Dual Formula]The dual of any finitely presented functor is representable.  More precisely, for any $F\in \fp(\A,\ab)$, $$F^*\cong (\blank,w(F))$$ where $F^*$ is the finitely presented functor in $\fp(\A^{op},\ab)$ given by $$F^*(A)=\nat(F,(A,\blank))$$ \end{cor*}  

In addition, we are able to completely describe the connection between the injective stabilization sequence and the defect sequence.  For any finitely presented functor $F$, there is a commutative diagram with exact rows \begin{center}\begin{tikzpicture}
\matrix(m)[ampersand replacement=\&, matrix of math nodes, row sep=3em, column sep=1.5em, text height=1.5ex, text depth=0.25ex]
{0\&\ext^1(\tr(F),\hom)\&F\&F^{**}\&\ext^2(\tr(F),\hom)\&0\\
0\&F_0\&F\&\big(w(F),\blank\big)\&F_1\&0\\
0\&\overline{F}\&F\&R^0F\&\tilde{F}\&0\\}; 
\path[->, font=\scriptsize]
(m-1-2) edge node[auto]{$\cong$}(m-2-2)
(m-2-2) edge node[auto]{$\cong$}(m-3-2)
(m-1-3) edge node[auto]{$1_F$}(m-2-3)
(m-2-3) edge node[auto]{$1_F$}(m-3-3)
(m-1-4) edge node[auto]{$\cong$}(m-2-4)
(m-2-4) edge node[auto]{$\cong$}(m-3-4)
(m-1-5) edge node[auto]{$\cong$}(m-2-5)
(m-2-5) edge node[auto]{$\cong$}(m-3-5);
\path[->, font=\scriptsize]
(m-1-1) edge node[auto]{}(m-1-2)
(m-1-2) edge node[auto]{}(m-1-3)
(m-1-3) edge node[auto]{}(m-1-4)
(m-1-4) edge node[auto]{}(m-1-5)
(m-1-5) edge node[auto]{}(m-1-6)
(m-2-1) edge node[auto]{}(m-2-2)
(m-2-2) edge node[auto]{}(m-2-3)
(m-2-3) edge node[auto]{$\varphi$}(m-2-4)
(m-2-4) edge node[auto]{}(m-2-5)
(m-2-5) edge node[auto]{}(m-2-6)
(m-3-1) edge node[auto]{}(m-3-2)
(m-3-2) edge node[auto]{}(m-3-3)
(m-3-3) edge node[auto]{}(m-3-4)
(m-3-4) edge node[auto]{}(m-3-5)
(m-3-5) edge node[auto]{}(m-3-6);
\end{tikzpicture}\end{center}  Thus the defect sequence, injective stabilization sequence, and the double dual sequence are all isomorphic.  From this we have the following formula.

\begin{thm*}[\fp-Injective Stabilization Formula]For any $F\in \fp(\A,\ab)$ 
\begin{enumerate}
\item $F_0\cong\overline{F}\cong\ext^1(\tr(F),\hom)$
\item $F_1\cong \tr\u\tr(F)\cong\ext^2(\tr(F),\hom)$ 
\item $\u\tr\u\tr F$ is the image of the unit of adjunction $\varphi=u_F\:F\to F^{**}$
\end{enumerate} \end{thm*}

By using Gentle's results from \cite{rongentle} concerning existence of injectives in $\fp(\A,\ab)$ whenever $\A$ has enough projectives and the \fp-dual formula, we are able to compute the derived functors $L^k(\blank)^*$.

\begin{thm*}The contravariant functor $L^0(\blank)^*$ is the only non-trivial left derived functor of the dual.  It is exact and its kernel is a Serre subcategory of $\fp(\Mod(R),\ab)$ which vanishes on functors whose presentations arise from certain short exact sequences. 
\end{thm*}

%
The structure of the paper is as follows.  In Section \ref{section2}, we review the Yoneda lemma and recall the definition of the category of finitely presented functors $\fp(\A,\ab)$ for an abelian category $\A$.  A functor $F\:\A\To \ab$ is  \dbf{finitely presented} if there exists a sequence of natural transformations $$(Y,\blank)\To (X,\blank)\To F\To 0$$  such that for any $A\in \A$ the sequence of abelian groups $$(Y,A)\To (X,A)\To F(A)\To 0$$ is exact.   The category $\fp(\A,\ab)$ is abelian and a sequence of functors $F\to G\to H$ is exact if and only if for every $A\in \A$, the sequence of abelian groups $F(A)\to G(A)\to H(A)$ is exact.

In Section \ref{section3}, we extend the definition of the injective stabilization of a functor $F\:\A\To\B$ first given by Auslander and Bridger in \cite{stab}.   Their definition requires that $\A$ is abelian and has enough injectives and that $\B$ is abelian.   We give a definition of injective stabilization of $F\:\A\To \B$ which does not require that $\A$ has enough injectives, but instead requires that a suitable zeroth right derived functor $R^0F$ exists.   Three definitions of zeroth right derived functors are discussed.  The first is the classical definition via injective resolutions and can be found in most introductory texts on homological algebra.  The second is due to Gabriel, appears in \cite{gab}, and requires that $\B$ is Grothendieck and has exact direct limits.  The third definition unifies the classical approach and Gabriel's approach by positing that the proper way to view the zeroth right derived functor of $F$ is as an approximation of $F$ by a left exact functor $R^0F$ arising from a particular unit of adjunction.  

Let $\S$ be a subcategory of a functor category $(\A,\B)$ where $\A,\B$ are abelian.  Denote by $\lex(\S)$ the full subcategory of $\S$ consisting of the left exact functors.  We say that $\S$ admits a \dbf{zeroth right derived functor} if the inclusion functor $s\:\lex(\S)\To \S$ has a left adjoint $r^0\:\S\to \lex(\S)$ such that 
\begin{enumerate}
\item The unit of adjunction $u\:1_{\S}\to sr^0$ is an isomorphism on the injectives of $\A$.  More precisely, if $F\in \S$ and $I\in \A$ is injective, then the morphism $(u_F)_I$ is an isomorphism.  
\item The composition $r^0s$ is isomorphic to the identity.  That is $$r^0s\cong 1_{\lex(\S)}$$
\end{enumerate} 
We close the section by using the more general definition of the zeroth right derived functor to define the injective stabilization sequence $$0\To \overline{F}\To F\To R^0F\To \tilde F\To 0$$ for any additive functor $F\:\A\To \B$  between abelian categories admitting a zeroth right derived functor.  This extends the definition of the injective stabilization sequence originally introduced by Auslander and Bridger in \cite{stab} and will allow us to compare the defect sequence with the injective stabilization sequence.

 In Section \ref{section4}, we recall Auslander's construction of the defect $$w\:\fp(\A,\ab)\To \ab$$ and discuss some of its properties.  Any finitely presented functor $F$ with presentation $$(Y,\blank)\To (X,\blank)\To F\To 0$$ gives rise to an exact sequence $$0\To w(F)\To X\To Y$$  On pages 202 through 205 in \cite{coh}, Auslander establishes that the assignment $F\mapsto w(F)$ determines an exact contravariant functor $w\:\fp(\A,\ab)\To \A$,  that the functor $w$ maps any representable $(X,\blank)$ to $X$, and that $w$ vanishes on those functors $F$ which arise from short exact sequences.  That is $w(F)=0$ if and only if there exists an exact sequence in $\A$ $$0\To X\To Y\To Z\To 0$$ such that we have an exact sequence $$0\To (Z,\blank)\To (Y,\blank)\To (X,\blank)\To F\To 0$$ We recall how Auslander constructs, for any finitely presented functor $F$, the defect sequence $$0\To F_0\To F\overset{\varphi}{\To} \big(w(F),\blank\big)\To F_1\To 0$$

The defect sequence plays an important role in establishing applications.  These are discussed in Section \ref{section5}. First we show that for any abelian category $\A$, the category of finitely presented functors $\fp(\A,\ab)$ admits a zeroth right derived functor as defined in Section \ref{section3} and in fact $$R^0F\cong \big(w(F),\blank\big)\cong F^{**}.$$ This result was the motivation behind introducing the more general definition of the zeroth right derived functor.  Auslander comments in \cite{coh} that if $\A$ has enough injectives, then $R^0F\cong \big(w(F),\blank\big)$. However, in the case that $\A$ does not have enough injectives, the defect sequence still exists and still induces an approximation of any finitely presented functor $F$ by the left exact functor $\big(w(F),\blank\big)$.  

Gabriel's definition of $R^0F$ is actually sufficient for showing that $R^0F\cong \big(w(F),\blank\big)$; however, by using the more general approach and focusing  on the map $F\overset{\varphi}\to \big(w(F),\blank\big)$ as an approximation, we are naturally led to  certain applications of the defect to the functor category.   The first of these is a direct consequence of Auslander's observation that for any finitely presented functor $F\:\A\To \ab$ and any left exact functor $G\:\A\To \ab$ there is an isomorphism $$\nat(F,G)\cong G\big(w(F)\big)$$ By setting $G=(X,\blank)$ we have the CoYoneda lemma.

From the CoYoneda Lemma we immediately get the \fp-dual formula which states that the dual of a finitely presented functor $F$ is $\big(\blank, w(F)\big)$ making it representable and hence projective.  Applying the \fp-dual formula to the tensor functor $X\otimes \blank$ whenever $X$ is finitely presented yields the following.  \begin{cor*}For any finitely presented module $X$, $$(X\otimes \blank\ )^*\cong(\blank, X^*)$$\end{cor*}

Finally, we show the \fp-injective stabilization formula which states that the injective stabilization of a finitely presented functor $F$ is $\ext^1(\tr(F),\hom)$.  In fact, we establish that for any abelian category $\A$ and for any finitely presented functor $F$, the injective stabilization sequence, defect sequence, and double dual sequence are all isomorphic.  In doing so we prove that $\u\tr\u\tr$ and $\tr\u\tr$ are endofunctors on the category $\fp(\A,\ab)$. This is surprising considering that these constructions are not in general functorial for the category of finitely presented modules, $\mod(R)$.

In Section \ref{section6}, we discuss the injective objects of $\fp(\A,\ab)$ when $\A$ has enough projectives.  In \cite{rongentle}, Gentle essentially shows that $\fp(\A,\ab)$ has enough injectives whenever $\A$ has enough projectives.   As a result both functor categories $\fp(\Mod(R),\ab)$ and $\fp(\Mod(R)^{op},\ab)$ have enough injectives.    We use the injectives of these categories to study the left derived functors of the dual $(\blank)^*$.  We show that the only non-trivial left derived functor of the dual is the functor $L^0(\blank)^*$ and that the pair of contravariant functors 
\begin{center}\begin{tikzpicture}
\matrix (m) [ampersand replacement= \&,matrix of math nodes, row sep=3em, 
column sep=6em,text height=1.5ex,text depth=0.25ex] 
{\fp(\Mod(R),\ab)\&\fp(\Mod(R)^{op},\ab)\\}; 
\path[->,thick, font=\scriptsize]
(m-1-1) edge[bend left] node[auto]{$L^0(\blank)^*$}(m-1-2)
(m-1-2) edge[bend left] node[below]{$L^0(\blank)^*$}(m-1-1);
\end{tikzpicture}\end{center}   are in fact exact.  The kernel of $L^0(\blank)^*$ is a Serre subcategory of $\fp(\Mod(R),\ab)$ and in fact, if $L^0(\blank)^*$ vanishes on $F$, then $w(F)=0$.    Therefore, the functor $L^0(\blank)^*$ detects some type of short exact sequence in $\Mod(R)$.  Using the defect and the \fp-dual formula, it is shown that $L^0(\blank)^*$ does not vanish exactly when $w$ vanishes.   Hence $L^0(\blank)^*$ detects a new type of short exact sequence.

In Section \ref{section7}, we demonstrate how almost split sequences are detected by  the defect.  The proof of the existence of these sequences is due to Auslander and Reiten and first appears in \cite{ar3}.  It turns out that almost split sequences can be predicted by looking at cases in which there exist minimal projective presentations of simple finitely presented functors $S$ for which $w(S)=0$.

\subsection*{Acknowledgements}

The author wishes to sincerely thank the reviewers for their many useful suggestions and comments concerning the original submission of the paper.  The analysis provided by the reviewers resulted in significant improvement to the structure of the paper and the presentation of the results.  In particular, the inclusion of the computations of $L^k(\blank)^*$ are directly due to the request by the reviewers to see applications of the results contained in the first submission.
%
%
%
\section{\label{section2}The Yoneda Lemma and Finitely Presented Functors}  
This section is a quick review.  For more details, see \cite{coh} and \cite{mclane}.
Throughout this entire paper:
\begin{enumerate}
\item The word functor will always mean additive functor.    
\item The category of right modules will be denoted by $\Mod(R)$ and the category of finitely presented right modules will be denoted by $\mod(R)$.  
\item Every left module can be viewed as a right module over its opposite ring $R^{op}$.
\item Unless otherwise stated, the category $\A$ will always be assumed to be an abelian category.  All statements hold for an arbitrary abelian category $\A$ which does not need to be skeletally small.
\item The category of abelian groups will be denoted by $\ab$.
\end{enumerate}

A functor $F\:\A\To \ab$ is called \dbf{representable} if it is isomorphic to $\hom_\A(X,\blank)$ for some $X\in \A$.   We will abbreviate the representable functors by $(X,\blank)$.   The most important property of representable functors is the following well known lemma of Yoneda:

\begin{lem}[Yoneda]For any covariant functor $F\:\A\to \ab$ and any $X\in \A$, there is an isomorphism:  $$\nat\big((X,\blank),F\big)\cong F(X)$$  given by $\alpha\mapsto \alpha_X(1_X)$.  The isomorphism is natural in both $F$ and $X$.  \end{lem}

An immediate consequence of the Yoneda lemma is that for any $X,Y\in \A$, $\nat\big((Y,\blank),(X,\blank)\big)\cong (X,Y)$.  Hence all natural transformations between representable functors come from morphisms between objects in $\A$.   Given a natural transformation between two functors $\alpha:F\To G$, there are functors $\coker(\alpha)$ and $\ker(\alpha)$ which are determined up to isomorphism by their value on any $A\in \A$ by the exact sequence in $\ab$:

$$0\To \ker(\alpha)(A)\To F(A)\overset{\alpha_A}{\To}G(A)\To \coker(\alpha)(A)\To 0$$

The following is a direct consequence of the Yoneda lemma and the fact that $\A$ is abelian.  For any morphism $\alpha\:(Y,\blank)\To (X,\blank)$, the functor $\ker(\alpha)$ is also representable.  Since $\alpha=(f,\blank)$ for some $f\:X\To Y$, applying  $(\blank,A)$ to the exact sequence \begin{center}\begin{tikzpicture}
\matrix(m)[ampersand replacement=\&, matrix of math nodes, row sep=3em, column sep=2.5em, text height=1.5ex, text depth=0.25ex]
{X\&Y\&Z\&0\\}; 
\path[->]
(m-1-1) edge node[auto]{$f$}(m-1-2)
(m-1-2) edge node[auto]{$g$}(m-1-3)
(m-1-3) edge node[auto]{}(m-1-4);
\end{tikzpicture}\end{center} results in the following exact sequence  \begin{center}\begin{tikzpicture}
\matrix(m)[ampersand replacement=\&, matrix of math nodes, row sep=3em, column sep=2.5em, text height=1.5ex, text depth=0.25ex]
{0\&(Z,A)\&(Y,A)\&(X,A)\\}; 
\path[->]
(m-1-1) edge node[auto]{}(m-1-2)
(m-1-2) edge node[auto]{$(g,A)$}(m-1-3)
(m-1-3) edge node[auto]{$(f,A)$}(m-1-4);
\end{tikzpicture}\end{center} Hence, the kernel of any natural transformation between representable functors is itself representable.  Specifically, $$\ker(\alpha)\cong (Z,\blank).$$

\begin{defn}A functor $F\:\A\To \ab$ is \dbf{finitely presented} if there exist $X,Y\in \A$ and $\alpha\:(Y,\blank)\To (X,\blank)$ such that $F\cong \coker(\alpha)$.\end{defn}

One easily shows that if $F$ is finitely presented, then the collection of natural transformations $\nat(F,G)$ for any functor $G\:\A\To \ab$ is actually an abelian group.  As such, one may form a category whose objects are the covariant finitely presented functors $F\:\A\To \ab$ and whose morphisms are the natural transformations between two such functors.  This category is denoted by $\fp(\A,\ab)$ and was studied extensively by Auslander in multiple works.  

\begin{thm}[Auslander, \cite{coh}, Theorem 2.3]The category $\fp(\A,\ab)$ consisting of all finitely presented functors together with the natural transformations between them satisfies the following properties:
\begin{enumerate}
\item $\fp(\A,\ab)$ is abelian.  A sequence of finitely presented functors $$F\To G\To H$$ is exact if and only if for every $A\in \A$ the sequence of abelian groups $$F(A)\To G(A)\To H(A)$$ is exact.

\item Every finitely presented functor $F\in \fp(\A,\ab)$ has a projective resolution of the form $$0\To (Z,\blank)\To (Y,\blank)\To (X,\blank)\To F\To 0$$
\end{enumerate}
\end{thm}

The functor $\Y\:\A\To \fp(\A,\ab)$ given by $\Y(X):=(X,\blank)$ is a contravariant left exact embedding commonly referred to as the \dbf{Yoneda}\dbf{ embedding}.  In addition, one can show that the projective objects of $\fp(\A,\ab)$ are exactly the representables.

\section{\label{section3}Zeroth Derived Functors and The Injective Stabilization}

Throughout this section we fix abelian categories $\A,\B$.  We will discuss three approaches to associating to each covariant functor $F\:\A\To \B$, a zeroth right derived functor $R^0F\:\A\To \B$.   The first is the classical construction which assumes that $\A$ has enough injectives.  The second comes from Gabriel's paper \cite{gab} on abelian categories and assumes that $\B$ has exact directed limits and is Grothendieck.  The final construction makes no assumptions on $\A$ or $\B$ other than that both are abelian. 

\subsection{Classical Definition}
Assume that $\A$ has enough injectives and let $F\:\A\To \ab$. The classical definition of $R^0F\:\A\To \B$ uses injective resolutions of the objects in $\A$.    For details, the reader is referred to \cite{ce}, Chapter 5, Section 5.   For any $X\in \A$, take any exact sequence $$0\to X\to I^0\to I^1$$ with $I^0,I^1$ injective.  The component $R^0F(X)$ is defined by the exact sequence

$$0\To R^0F(X)\To F(I^0)\To F(I^1) $$

It is easily seen that this assignment is functorial in both $F$ and $X$.  Moreover, up to isomorphism $R^0F$ is independent of the choices of $I^0$ and $I^1$.   The following properties may also be established and the reader is again referred to \cite{ce}: 
\begin{enumerate}
\item $R^0F$ and $F$ agree on injectives.
\item $R^0F$ is left exact.
\item If $F$ is left exact, then $R^0F\cong F$.
\item If $G$ is left exact and $G$ agrees on injectives with $F$, then $G\cong R^0F$.
\item There is a natural transformation $\eta\:F\to R^0F$ satisfying the following universal property.  For any  $\alpha\: F\to G$ with $G$ left exact, there exists a unique  $\phi\:R^0F\to G$ making that the  following diagram commute
\dia{F\ar[r]^{\eta\quad}\ar[d]_{\forall\alpha}&R^0F\ar[dl]^{\exists !\phi}\\
G}

\item For any left exact functor $G$, there exists an isomorphism $$(F,G)\cong (R^0F,G)$$ which is natural in $F$ and $G$.
\end{enumerate}

\subsection{Gabriel's Definition}In \cite{gab}, Gabriel constructs $R^0F$ under the assumption that the source category is abelian and that the target category is Grothendieck and has exact directed limits.  The following is the construction given there. Let $A\in \A$ and define $\mathbb{E}_A$ to be the set of all short exact sequences $0\to A\to B\to C\to 0$ where $B,C\in \A$.    Define the relation $\le$ on $\mathbb{E}_A$ by $$0\to A\to B\to C\to 0\qquad \le \qquad 0\to A\to B'\to C'\to 0$$ if and only if there exists commutative diagram
\dia{0\ar[r]&A\ar@{=}[d]\ar[r]&B\ar[d]\ar[r]&C\ar[d]\ar[r]&0\\
0\ar[r]&A\ar[r]& B'\ar[r]&C'\ar[r]&0}
For any two short exact sequences $\eps,\eps'\in\mathbb{E}_A$: 
\dia{\eps:0\ar[r] &A\ar[r]^i& B\ar[r]& C\ar[r]& 0} 
\dia{\eps':0\ar[r] & A\ar[r]^j &B'\ar[r]&C'\ar[r]& 0} the pushout \dia{A\ar[d]_j\ar[r]^i&B\ar[d]^b\\
B'\ar[r]_{b'}&E} gives rise to the short exact sequence $\eps''\in\mathbb{E}_A$ \dia{\eps'':0\ar[r]&A\ar[r]^{bi=b'j}&E\ar[r]&U\ar[r]&0}   satisfying $\eps,\eps'\le \eps''$. 

Let \dia{0\ar[r]&A\ar[r]&B\ar[r]^f&C\ar[r]&0} be exact in $\A$.  For any functor $F\:\A\To\B$, there is an exact sequence \dia{0\ar[r]&K\ar[r]&F(B)\ar[r]^{F(f)}&F(C)}   Denote by $F(\mathbb{E}_A)$  all exact sequences in $\B$ arising in this way.  The relation $\le$ on $\mathbb{E}_A$ induces a relation $\le_F$ on $F(\mathbb{E}_A)$ given by

$$0\To K\To F(B)\overset{F(f)}\To F(C)\qquad \le_F \qquad 0\To K'\To F(B')\overset{F(f')}\To F(C')$$ if and only if there exists a commutative diagram

\dia{0\ar[r]&K\ar[d]\ar[r]&F(B)\ar[d]\ar[r]^{F(f)}&F(C)\ar[d]\ar[r]&0\\
0\ar[r]&K'\ar[r]& F(B')\ar[r]^{F(f')}&F(C')\ar[r]&0}  The relation $\le_F$ makes $F(\mathbb{E}_A)$ into a directed system.  Gabriel's definition of $R^0F$ at $A$ is $$R^0F(A)=\underset{\To}{\lim}\ F(\mathbb{E}_A)$$  This completely determines $R^0F$ as a functor up to isomorphism.  

\begin{thm}[Gabriel, \cite{gab}]For any $F\:\A\To \B$,
\begin{enumerate}
\item $R^0F$ and $F$ agree on injectives.
\item $R^0F$ is left exact.
\item If $F$ is left exact, then $R^0F\cong F$.
\item There is a natural transformation $\eta\:F\to R^0F$ satisfying the following universal property.  For any  $\alpha\: F\to G$ with $G$ left exact, there exists a unique  $\phi\:R^0F\to G$ making that the  following diagram commute
\dia{F\ar[r]^{\eta\quad}\ar[d]_{\forall\alpha}&R^0F\ar[dl]^{\exists !\phi}\\
G}

\item For any left exact functor $G$, there exists an isomorphism $$(F,G)\cong (R^0F,G)$$ which is natural in $F$ and $G$.
\end{enumerate}
\end{thm}

\subsection{More General Definition}
Let $\S$ denote any subcategory of $(\A,\B)$.  Define $\lex(\S)$ to be the full subcategory of $\S$ consisting of the left exact functors, that is, all functors $F\in \S$ such that if $0\to A\to B\to C\to 0$ is exact, then $0\to F(A)\to F(B)\to F(C)$ is exact.  There is an inclusion functor  \begin{center}\begin{tikzpicture}
\matrix(m)[ampersand replacement=\&, matrix of math nodes, row sep=3em, column sep=2.5em, text height=1.5ex, text depth=0.25ex]
{\lex(\S)\&\S\\}; 
\path[->]
(m-1-1) edge node[auto]{$s$}(m-1-2);
\end{tikzpicture}\end{center}  

\begin{defn}We say that $\S$ admits a \dbf{zeroth right derived functor} if $s$ has a left adjoint $r^0\:\S\to \lex(\S)$ such that 
\begin{enumerate}
\item The unit of adjunction $u\:1_{\S}\to sr^0$ is an isomorphism on the injectives of $\A$.  More precisely, if $F\in \S$ and $I\in \A$ is injective, then the morphism $(u_F)_I$ is an isomorphism.  
\item The composition $r^0s$ is isomorphic to the identity.  That is $$r^0s\cong 1_{\lex(\S)}$$
\end{enumerate}\end{defn}

If the inclusion $s$ has a left adjoint, the category $\lex(\S)$ is a reflective subcategory of $\S$.   If $\S$ admits a zeroth right derived functor, then we define the composition $R^0=sr^0$ to be the \dbf{zeroth right} \dbf{derived functor of $\S$}.  Clearly if $\S\sub (\A,\B)$ admits a zeroth right derived functor $R^0$,  then for any functor $F\in \S$ and any injective $I\in \A$, $$R^0F(I)\cong F(I)$$  Moreover, the functor $r^0\:\S\to \lex(\S)$ produces for each functor $F$ a left exact functor $r^0F$ by altering $F$ in the smallest amount possible. To clarify this comment, note that $r^0$ does not change the functors which are already left exact.  It also does not change values of $F$ on objects $X$ such that $0\to X\to Y\to Z\to 0$ splits for all $Y,Z$ as this condition implies that $0\to F(X)\to F(Y)\to F(Z)$ is exact whether or not $F$ is left exact.  These objects are precisely the injectives of $\A$.  If $\A$ has enough injectives,  it is easily seen that the classical definition of $R^0$ satisfies the more general definition, that is $(\A,\B)$ admits a zeroth right derived functor.   Similarly, if $\B$ is Grothendieck and has exact directed limits, Gabriel's definition also satisfies the more general definition.

\subsection{The Injective Stabilization Sequence}
Now suppose that $\S$ is an abelian subcategory of $(\A,\B)$ that admits a zeroth right derived functor $R^0\:\S\to \S$.  In this case, there is an exact sequence of functors  \begin{center}\begin{tikzpicture}
\matrix(m)[ampersand replacement=\&, matrix of math nodes, row sep=3em, column sep=2.5em, text height=1.5ex, text depth=0.25ex]
{0\&\ker(u)\&1_{\S}\&R^0\&\coker(u)\&0\\}; 
\path[->]
(m-1-1) edge node[auto]{}(m-1-2)
(m-1-2) edge node[auto]{}(m-1-3)
(m-1-3) edge node[auto]{$u$}(m-1-4)
(m-1-4) edge node[auto]{}(m-1-5)
(m-1-5) edge node[auto]{}(m-1-6);
\end{tikzpicture}\end{center}where $u$ is the unit of adjunction.  The functor $\ker(u)$ is called the \dbf{injective } \dbf{stabilization functor}.  Given $F\in \S$, the functor $\ker(u)(F)$ will be denoted $\overline{F}$. The functor $\overline{F}$ is called the \dbf{injective }\dbf{stabilization of F}.  Evaluating the exact sequence at $F$ yields an exact sequence of functors called the \dbf{injective stabilization sequence}:\begin{center}\begin{tikzpicture}
\matrix(m)[ampersand replacement=\&, matrix of math nodes, row sep=3em, column sep=2.5em, text height=1.5ex, text depth=0.25ex]
{0\&\overline{F}\&F\&R^0F\&\tilde F\&0\\}; 
\path[->]
(m-1-1) edge node[auto]{}(m-1-2)
(m-1-2) edge node[auto]{}(m-1-3)
(m-1-3) edge node[auto]{$u_F$}(m-1-4)
(m-1-4) edge node[auto]{}(m-1-5)
(m-1-5) edge node[auto]{}(m-1-6);
\end{tikzpicture}\end{center}  A functor $F\in \S$ is called \dbf{injectively stable} if $R^0F=0$.   This generalizes the definition of the injective stabilization given by Auslander and Bridger in \cite{stab}.  The definition given there is essentially the same except that it uses the zeroth right derived functor as defined classically, which requires the existence of injectives in $\A$.

\section{\label{section4}The Defect of a Finitely Presented Functor}

 At this point we will explicitly recall the construction of the \dbf{defect} functor $w\:\fp(\A,\ab)\to \A$, originally introduced by Auslander in \cite{coh}.  Let $F\in \fp(\A,\ab)$ and let  
\begin{center}\begin{tikzpicture}
\matrix(m)[ampersand replacement=\&, matrix of math nodes, row sep=3em, column sep=2.5em, text height=1.5ex, text depth=0.25ex]
{(Y,\blank)\&(X,\blank)\&F\&0\\}; 
\path[->, font=\scriptsize]
(m-1-1) edge node[auto]{$(f,\blank)$}(m-1-2)
(m-1-2) edge node[auto]{$\alpha$}(m-1-3)
(m-1-3) edge node[auto]{}(m-1-4);
\end{tikzpicture}\end{center} be any presentation. The value of $w$ at $F$ is defined by the exact sequence 
\begin{center}\begin{tikzpicture}
\matrix(m)[ampersand replacement=\&, matrix of math nodes, row sep=3em, column sep=2.5em, text height=1.5ex, text depth=0.25ex]
{0\&w(F)\&X\&Y\\}; 
\path[->, font=\scriptsize]
(m-1-1) edge node[auto]{}(m-1-2)
(m-1-2) edge node[auto]{$k$}(m-1-3)
(m-1-3) edge node[auto]{$f$}(m-1-4);
\end{tikzpicture}\end{center}This assignment extends to a contravariant additive functor $$w\:\fp(\A,\ab)\to \A$$  Auslander established the following properties concerning the functor $w$:
\begin{enumerate}
\item For any finitely presented functor $F$, $w(F)$ is independent of the chosen projective presentation.
\item $w(X,\blank)\cong X$.
\item For any finitely presented functor $F$, $w(F)=0$ if and only if there exists an exact sequence $0\to X\to Y\to Z\to 0$ such that $0\to(Z,\blank)\to (Y,\blank)\to (X,\blank)\to F\to 0$ is exact.
\item The functor $w$ is exact.\end{enumerate}

We now recall the construction of the \dbf{defect sequence}  \begin{center}\begin{tikzpicture}
\matrix(m)[ampersand replacement=\&, matrix of math nodes, row sep=3em, column sep=2.5em, text height=1.5ex, text depth=0.25ex]
{0\&F_0\&F\&\big(w(F),\blank\big)\&F_1\&0\\}; 
\path[->, font=\scriptsize]
(m-1-1) edge node[auto]{}(m-1-2)
(m-1-2) edge node[auto]{}(m-1-3)
(m-1-3) edge node[auto]{$\varphi$}(m-1-4)
(m-1-4) edge node[auto]{}(m-1-5)
(m-1-5) edge node[auto]{}(m-1-6);
\end{tikzpicture}\end{center}  To any projective resolution  \begin{center}\begin{tikzpicture}
\matrix(m)[ampersand replacement=\&, matrix of math nodes, row sep=3em, column sep=2.5em, text height=1.5ex, text depth=0.25ex]
{0\&(Z,\blank)\&(Y,\blank)\&(X,\blank)\&F\&0\\}; 
\path[->, font=\scriptsize]
(m-1-1) edge node[auto]{}(m-1-2)
(m-1-2) edge node[auto]{$(g,\blank)$}(m-1-3)
(m-1-3) edge node[auto]{$(f,\blank)$}(m-1-4)
(m-1-4) edge node[auto]{$\alpha$}(m-1-5)
(m-1-5) edge node[auto]{}(m-1-6);
\end{tikzpicture}\end{center} apply the exact functor $w$.  This yields the exact sequence \begin{center}\begin{tikzpicture}
\matrix(m)[ampersand replacement=\&, matrix of math nodes, row sep=3em, column sep=2.5em, text height=1.5ex, text depth=0.25ex]
{0\&w(F)\&X\&Y\&Z\&0\\}; 
\path[->, font=\scriptsize]
(m-1-1) edge node[auto]{}(m-1-2)
(m-1-2) edge node[auto]{$k$}(m-1-3)
(m-1-3) edge node[auto]{$f$}(m-1-4)
(m-1-4) edge node[auto]{$g$}(m-1-5)
(m-1-5) edge node[auto]{}(m-1-6);
\end{tikzpicture}\end{center}  This exact sequence embeds into the following commutative diagram with exact rows and columns: 
\dia{&&0\ar[d]\\
&w(F)\ar@{=}[d]\ar@{=}[r]&w(F)\ar[d]^k\\
0\ar[r]&w(F)\ar[r]^k&X\ar[d]\ar[r]&Y\ar@{=}[d]\ar[r]&Z\ar@{=}[d]\ar[r]&0\\
&0\ar[r]&V\ar[d]\ar[r]&Y\ar[r]&Z\ar[r]&0\\
&&0} Applying the Yoneda embedding to this diagram and extending to include cokernels where necessary yields the following commutative diagram with exact rows and columns: \dia{&&&0\ar[d]&0\ar[d]\\
0\ar[r]&(Z,\blank)\ar@{=}[d]\ar[r]&\ar@{=}[d](Y,\blank)\ar[r]&(V,\blank)\ar[d]\ar[r]&F_0\ar[d]\ar[r]&0\\
0\ar[r]&(Z,\blank)\ar[r]&(Y,\blank)\ar[r]&(X,\blank)\ar[d]_{(k,\blank)}\ar[r]^{\alpha}&F\ar[d]^{\varphi}\ar[r]&0\\
&&&(w(F),\blank)\ar[d]\ar@{=}[r]&(w(F),\blank)\ar[d]\\
&&&F_1\ar[d]\ar@{=}[r]&F_1\ar[d]\\
&&&0&0}  This yields the following exact sequence: 
\begin{center}\begin{tikzpicture}
\matrix(m)[ampersand replacement=\&, matrix of math nodes, row sep=3em, column sep=2.5em, text height=1.5ex, text depth=0.25ex]
{0\&F_0\&F\&\big(w(F),\blank\big)\&F_1\&0\\}; 
\path[->, font=\scriptsize]
(m-1-1) edge node[auto]{}(m-1-2)
(m-1-2) edge node[auto]{}(m-1-3)
(m-1-3) edge node[auto]{$\varphi$}(m-1-4)
(m-1-4) edge node[auto]{}(m-1-5)
(m-1-5) edge node[auto]{}(m-1-6);
\end{tikzpicture}\end{center} where $w(F_0)=0$ and $w(F_1)=0$.  This sequence is functorial in $F$. The map $\varphi$ is the unique map such that $\varphi \alpha=(k,\blank)$.

Let $\fp_0(\A,\ab)$ denote the full subcategory of $\fp(\A,\ab)$ consisting of all functors $F$ for which $w(F)=0$.  It is easily seen that $\fp_0(\A,\ab)$ is an abelian category and that the embedding $$\fp_0(\A,\ab)\To \fp(\A,\ab)$$ is exact and reflects exact sequences.  Moreover, $F\in\fp_0(\A,\ab)$ if and only if $F_0\cong F$.  As a result, for any $F\in \fp_0(\A,\ab)$ and for any $G\in \fp(\A,\ab)$, there exists an isomorphism $$\nat(F,G)\cong \nat(F,G_0)$$ natural in both $F$ and $G$.  This can stated more precisely as follows.

\begin{prop}The functor $F\mapsto F_0$ is the right adjoint to the exact embedding $$\fp_0(\A,\ab)\To \fp(\A,\ab)$$  As a right adjoint, this functor is left exact. Since this embedding is exact the endofunctor $F\mapsto F_0$ is also left exact.  \end{prop}

\section{\label{section5}The CoYoneda Lemma and Applications}
We begin this section by making the observation that if $F\in\fp(\A,\ab)$ and $w(F)=0$, then $F$ vanishes on injectives.  To see this note that since $w(F)=0$ there exists an exact sequence\begin{center}\begin{tikzpicture}
\matrix(m)[ampersand replacement=\&, matrix of math nodes, row sep=3em, column sep=2.5em, text height=1.5ex, text depth=0.25ex]
{0\&X\&Y\&Z\&0\\}; 
\path[->, font=\scriptsize]
(m-1-1) edge node[auto]{}(m-1-2)
(m-1-2) edge node[auto]{$f$}(m-1-3)
(m-1-3) edge node[auto]{$g$}(m-1-4)
(m-1-4) edge node[auto]{}(m-1-5);
\end{tikzpicture}\end{center} such that the following sequence is exact \begin{center}\begin{tikzpicture}
\matrix(m)[ampersand replacement=\&, matrix of math nodes, row sep=3em, column sep=2.5em, text height=1.5ex, text depth=0.25ex]
{0\&(Z,\blank)\&(Y,\blank)\&(X,\blank)\&F\&0\\}; 
\path[->, font=\scriptsize]
(m-1-1) edge node[auto]{}(m-1-2)
(m-1-2) edge node[auto]{$(g,\blank)$}(m-1-3)
(m-1-3) edge node[auto]{$(f,\blank)$}(m-1-4)
(m-1-4) edge node[auto]{$\alpha$}(m-1-5)
(m-1-5) edge node[auto]{}(m-1-6);
\end{tikzpicture}\end{center}  By evaluating this sequence on an injective we have the following exact sequence
\begin{center}\begin{tikzpicture}
\matrix(m)[ampersand replacement=\&, matrix of math nodes, row sep=3em, column sep=2.5em, text height=1.5ex, text depth=0.25ex]
{0\&(Z,I)\&(Y,I)\&(X,I)\&F(I)\&0\\}; 
\path[->, font=\scriptsize]
(m-1-1) edge node[auto]{}(m-1-2)
(m-1-2) edge node[auto]{$(g,I)$}(m-1-3)
(m-1-3) edge node[auto]{$(f,I)$}(m-1-4)
(m-1-4) edge node[auto]{$\alpha_I$}(m-1-5)
(m-1-5) edge node[auto]{}(m-1-6);
\end{tikzpicture}\end{center}   Since $I$ is injective, $(f,I)$ is an epimorphism.  Therefore $\alpha_I=0$.  Since $\alpha_I$ is an epimorphism, $F(I)=0$.

On page 204 of \cite{coh}, Auslander observes that the morphism $F\overset{\varphi}\To \big(w(F),\blank\big)$ satisfies the following universal property.  Given any left exact functor $G\:\A\To \ab$ and any natural transformation $\alpha\: F\To G$, there exists a natural transformation morphism $\phi\:\big(w(F),\blank\big)\To G$ such that the following diagram commutes

\dia{F\ar[r]^{\varphi\qquad}\ar[d]_{\forall\alpha}&\big(w(F),\blank\big)\ar[dl]^{\exists !\phi}\\
G}

From this observation and the Yoneda lemma one easily has the following.   
\begin{thm}[Auslander, \cite{coh}, page 204]\label{ref248}Let $F\in \fp(\A,\ab)$.  Suppose that $G\: \A\To\ab$ is  any left exact additive covariant functor.  Then there is an isomorphism $$(F,G)\cong G\big(w(F)\big)$$ which is natural in $F$ and $G$. \end{thm}

\begin{thm}The inclusion functor $s\:\lex(\fp(\A,\ab))\To \fp(\A,\ab)$ admits a left adjoint $r^0\:\fp(\A,\ab)\to \lex(\fp(\A,\ab))$.  Moreover, the unit of adjunction evaluated at any finitely presented functor $F$ is the map $\varphi$ in the defect sequence \begin{center}\begin{tikzpicture}
\matrix(m)[ampersand replacement=\&, matrix of math nodes, row sep=3em, column sep=2.5em, text height=1.5ex, text depth=0.25ex]
{0\&F_0\&F\&\big(w(F),\blank\big)\&F_1\&0\\}; 
\path[->, font=\scriptsize]
(m-1-1) edge node[auto]{}(m-1-2)
(m-1-2) edge node[auto]{}(m-1-3)
(m-1-3) edge node[auto]{$\varphi$}(m-1-4)
(m-1-4) edge node[auto]{}(m-1-5)
(m-1-5) edge node[auto]{}(m-1-6);
\end{tikzpicture}\end{center} \end{thm}

\begin{pf}Define $r^0(F):=\big(w(F),\blank\big)$.  From the preceding theorem, for any $F\in \fp(\A,\ab)$ and for any $G\in \lex(\fp(\A,\ab))$, there exists an isomorphism $$\big(F,s(G)\big)\cong (r^0F,G)$$  Moreover, this isomorphism is natural in $F$ and $G$  which establishes the adjunction. 

In the commutative diagram with exact rows: \begin{center}\begin{tikzpicture}
\matrix(m)[ampersand replacement=\&, matrix of math nodes, row sep=3em, column sep=4em, text height=1.5ex, text depth=0.25ex]
{0\&G\big(w(F)\big)\&G(X)\&G(Y)\\
0\&\big((w(F),\blank),G\big)\&\big((X,\blank),G\big)\&\big((Y,\blank),G\big)\\
0\&(F,G)\&\big((X,\blank),G\big)\&\big((Y,\blank),G\big)\\}; 
\path[->, font=\scriptsize]
(m-1-1) edge node[auto]{}(m-1-2)
(m-1-2) edge node[auto]{$G(k)$}(m-1-3)
(m-1-3) edge node[auto]{$G(f)$}(m-1-4)
(m-2-1) edge node[auto]{} (m-2-2)
(m-2-2) edge node[auto]{$\big((k,\blank),G\big)$} (m-2-3)
(m-2-3) edge node[auto]{$\big((f,\blank),G\big)$} (m-2-4)
(m-1-2) edge node[auto]{$\cong$} (m-2-2)
(m-1-3) edge node[auto]{$\cong$} (m-2-3)
(m-1-4) edge node[auto]{$\cong$} (m-2-4)
(m-3-1) edge node[auto]{} (m-3-2)
(m-3-2) edge node[auto]{$(\alpha,G)$} (m-3-3)
(m-3-3) edge node[auto]{$\big((f,\blank),G\big)$} (m-3-4)
(m-2-3) edge node[auto]{$1$}(m-3-3)
(m-2-4) edge node[auto]{$1$}(m-3-4)
(m-2-2) edge node[auto]{$\theta_{F,G}$} (m-3-2);
\end{tikzpicture}\end{center}  choose $G=\big(w(F),\blank\big)$ and focus on the lower left commutative square:  

\begin{center}\begin{tikzpicture}
\matrix(m)[ampersand replacement=\&, matrix of math nodes, row sep=3em, column sep=8em, text height=1.5ex, text depth=0.25ex]
{\bigg(\big(w(F),\blank\big),\big(w(F),\blank\big)\bigg)\&\bigg(\big(X,\blank\big),\big(w(F),\blank\big)\bigg)\\
\bigg(F,\big(w(F),\blank\big)\bigg)\&\bigg(\big(X,\blank\big),\big(w(F),\blank\big)\bigg)\\}; 
\path[->, font=\scriptsize]
(m-1-1) edge node[auto]{$\bigg((k,\blank),\big(w(F),\blank\big)\bigg)$} (m-1-2)
(m-1-1) edge node[left]{$\theta$} (m-2-1)
(m-1-2) edge node[right]{$1$} (m-2-2)
(m-2-1) edge node[below]{$\bigg(\alpha,\big(w(F),\blank\big)\bigg)$}(m-2-2);
\end{tikzpicture}\end{center}  By definition, the unit of adjunction evaluated at $F$ is $u_F=\theta(1)$. From this commutative square, it follows that $u_F\alpha=(k,\blank)$.  Since $\varphi$ is the unique map such that $\varphi\alpha=(k,\blank)$, it follows that $\varphi=u_F$.    $\qed$\end{pf}

\begin{thm}The category $\fp(\A,\ab)$ admits a zeroth right derived functor  $R^0\:\fp(\A,\ab)\to \fp(\A,\ab)$ and $R^0\cong \Y w$.  In particular the zeroth right derived functor applied to a finitely presented functor yields a representable functor.  \end{thm}

\begin{pf}Let $s\:\lex(\fp(\A,\ab))\To \fp(\A,\ab)$ be the natural inclusion functor.  Observe that for any functor $F$, $\Y w(F)=(w(F),\blank)$.  We have already shown that $s$ admits a left adjoint $r^0\:\fp(\A,\ab)\to \lex(\fp(\A,\ab))$ sending $F$ to $\big(w(F),\blank\big)$ and that the unit of this adjunction evaluated at $F$ is the map $\varphi$ in the exact sequence \begin{center}\begin{tikzpicture}
\matrix(m)[ampersand replacement=\&, matrix of math nodes, row sep=3em, column sep=2.5em, text height=1.5ex, text depth=0.25ex]
{0\&F_0\&F\&\big(w(F),\blank\big)\&F_1\&0\\}; 
\path[->, font=\scriptsize]
(m-1-1) edge node[auto]{}(m-1-2)
(m-1-2) edge node[auto]{}(m-1-3)
(m-1-3) edge node[auto]{$\varphi$}(m-1-4)
(m-1-4) edge node[auto]{}(m-1-5)
(m-1-5) edge node[auto]{}(m-1-6);
\end{tikzpicture}\end{center}  Since $w(F_0)=w(F_1)=0$, both $F_0$ and $F_1$ vanish on injectives.  Therefore $\varphi_I$ is an isomorphism whenever $I$ is injective. Finally, suppose that $F$ is left exact.  Then $F$ is easily seen to be representable.  Observe that if $F\cong (X,\blank)$, then $w(F)\cong X$.  Therefore $$r^0F\cong \big(w(F),\blank\big)\cong(X,\blank)\cong F\qquad$$ these isomorphisms being natural.  As a result, $r^0(F)\cong F$.  $\qed$
\end{pf}

The following is a direct consequence of Theorem \ref{ref248}.  It is clear from \cite{coh} that Auslander was aware of the result.  Krause states this result in \cite{kraus} in terms of the Yoneda embedding admitting an adjoint.

\begin{lem}[The CoYoneda Lemma]For any $F\in \fp(\A,\ab)$ and any $X\in \A$, $$\nat\big(F,(X,\blank)\big)\cong \big(X,w(F)\big)$$ this isomorphism being natural in $F$ and $X$.  This makes $(w,\Y)$ an adjoint pair.\end{lem}

\begin{pf}Since $(X,\blank)$ is left exact finitely presented, by Theorem \ref{ref248}, $$\nat\big(F,(X,\blank)\big)\cong (X,\blank)(w(F))=\big(X,w(F)\big)$$ this isomorphism being natural in $F$ and $X$. $\qed$\end{pf}

One of the advantages of studying functors is that one may view them as generalizations of modules, or alternatively, one may view modules as examples of functors.  Many of the familiar module-theoretic definitions can be translated to the functor category.  One such example is the analog of the functor $(\blank, R)=(\blank)^*$ which sends left modules to right modules and vice versa.

In \cite{fisheradjoint}, Fisher-Palmquist and Newell define a functor $$(\blank)^*\:(\A,\ab)\to (\A^{op},\ab)$$ as follows:  For each $F\:\A\to \ab$, $F^*\:\A^{op}\to \ab$ is defined by $$F^*(X):=\nat(F,(X,\blank))$$ It is easily seen that for any representable functor $(A,\blank)$, $$(A,\blank)^*=(\blank, A)$$  As a corollary to the CoYoneda lemma, the dual of every finitely presented functor will be representable and completely determined by the defect. At this point, it is unknown if the condition that $F^*$ is representable is necessary for $F$ to be finitely presented. 

\begin{cor}[\fp-Dual Formula]For any finitely presented functor $F$, $$F^*\cong (\blank,w(F)).$$\end{cor}

\begin{pf}Let $F\in \fp(\A,\ab)$. Then by definition $$F^*(X)\cong\nat(F,(X,\blank))$$ and by the CoYoneda lemma, $$\nat(F,(X,\blank))\cong (X,w(F)).$$   It follows that $F^*(X)\cong (X,w(F))$ on component $X$.  Therefore $$F^*\cong \big(\blank, w(F)\big)\qquad\qed$$\end{pf}

\begin{cor}For the category $\fp(\A,\ab)$, $R^0=(\blank)^{**}$\end{cor}

\begin{pf} $$F^{**}\cong\big(\blank,w(F)\big)^*\cong \big(w(F),\blank\big)\qquad\qed$$   \end{pf}

\begin{cor}If $F\in \fp(\A,\ab)$, then $F^*=0$ if and only if $w(F)=0$.\end{cor}

Let $X$ be a finitely presented right $R$-module.  From Lemma 6.1 in \cite{coh}, the functor $X\otimes \blank$ is finitely presented and if  $$Q\To P\To X\To 0$$ is a presentation of $X$ by finitely generated projectives $P,Q$, then $$(Q^*,\blank)\To (P^*,\blank)\To X\otimes \blank\To 0$$ is a presentation of the functor $X\otimes \blank$.  From this one can easily conclude that $w(\blank\otimes X)\cong X^*$ as follows.   Applying the exact functor $w$ to the presentation of $X\otimes \blank$ yields the following commutative diagram with exact rows $$0\To w(X\otimes \blank)\To P^*\To Q^*$$ Since $(\blank, R)$ is left exact we also have an exact sequence $$0\To X^*\To P^*\To Q^*$$  It  follows that $w(X\otimes \blank)\cong X^*$.    From the \fp-dual formula, we now conclude the following formula for the dual of a finitely presented tensor functor.

\begin{cor}\label{tdf}For any finitely presented module $X$, $$(X\otimes\blank\ )^*\cong (\blank, X^*).$$\end{cor}

To state our next result we will need the notion of the transpose of a finitely presented module.  For any ring $R$ and any module $M\in \mod(R)$, apply the functor $(\blank)^*=(\blank,R)$ to any presentation of $M$:
\begin{center}\begin{tikzpicture}
\matrix(m)[ampersand replacement=\&, matrix of math nodes, row sep=3em, column sep=2.5em, text height=1.5ex, text depth=0.25ex]
{P_1\&P_0\&X\&0\\};
\path[->]
(m-1-1) edge (m-1-2)
(m-1-2) edge (m-1-3)
(m-1-3) edge (m-1-4);
\end{tikzpicture}\end{center} yielding exact sequence
\begin{center}\begin{tikzpicture}
\matrix(m)[ampersand replacement=\&, matrix of math nodes, row sep=3em, column sep=2.5em, text height=1.5ex, text depth=0.25ex]
{0\&X^*\&(P_0)^*\&(P_1)^*\&\tr(M)\&0\\};
\path[->]
(m-1-1) edge (m-1-2)
(m-1-2) edge (m-1-3)
(m-1-3) edge (m-1-4)
(m-1-4) edge (m-1-5)
(m-1-5) edge (m-1-6);
\end{tikzpicture}\end{center} The module $\tr(M)$ is called the transpose of $M$.  The assignment $M\mapsto \tr(M)$ is not in general a functor and depends on the choices involved.  The construction first appears in \cite{coh}, though the notation $\tr(M)$ is not used there.     

\begin{prop}[Auslander, \cite{coh}, Proposition 6.3]\label{dseq}Let $R$ be a ring.  For any finitely presented module $M$, there is an exact sequence of functors $$0\to \ext^1(\tr(M),R)\to M\to M^{**}\to \ext^2(\tr(M),R)\to 0$$\end{prop}

We will refer to the sequence in Proposition \ref{dseq} as the \dbf{double dual sequence}.  It will have a functor analog in the category $\fp(\A,\ab)$.   First we recall the definition of the transpose of a finitely presented functor.  Given $F\in \fp(\A,\ab)$, we take presentation of $F$  

\begin{center}\begin{tikzpicture}
\matrix(m)[ampersand replacement=\&, matrix of math nodes, row sep=3em, column sep=2.5em, text height=1.5ex, text depth=0.25ex]
{0\&(Z,\blank)\&(Y,\blank)\&(X,\blank)\&F\&0\\}; 
\path[->, font=\scriptsize]
(m-1-1) edge node[auto]{}(m-1-2)
(m-1-2) edge node[auto]{}(m-1-3)
(m-1-3) edge node[auto]{}(m-1-4)
(m-1-4) edge node[auto]{}(m-1-5)
(m-1-5) edge node[auto]{}(m-1-6);
\end{tikzpicture}\end{center} Apply the left exact functor $(\blank)^*$ to get the exact sequence $$0\To F^*\To (\blank,X)\To (\blank,Y)\To \tr( F)\To 0$$  The functor $\tr(F)$ is not uniquely determined as different choices of presentations may yield different functors; however, it is easily shown that any two such functors are projectively equivalent.

For any ring $R$ and any right module $M\in \Mod(R)$ there is a natural transformation of functors $$M\otimes \blank\To (M^*,\blank)$$  Evaluating these functors at the ring $R$ gives a natural morphism of right modules $$M\To M^{**}$$

Let $F\:\A\to \ab$ and $G\:\A^{op}\to \ab$ be two functors.  There is a well known abelian group $F\otimes G$ which is a tensor product of functors analogous to the tensor product of modules.  This abelian group was studied extensively in both \cite{fishertensor} and \cite{fisheradjoint}.  

\begin{prop}[Fisher, \cite{fishertensor}]The tensor product $F\otimes G$ of functors is completely determined by the following properties:
\begin{enumerate}
\item For any $A\in \A$, $(A,\blank)\otimes G\cong G(A)$
\item For any $B\in \A$, $F\otimes(\blank,B)\cong F(B)$
\item For any $F\in (\A,\ab)$, the functor $F\otimes \blank$ is cocontinuous.
\item For any $G\in (\A^{op},\ab)$, the functor $\blank\otimes G$ is cocontinuous.
\end{enumerate}
\end{prop}

The natural morphism from a module to its double dual is obtained by evaluating the natural morphism $M\otimes\blank \To (M^*,\blank)$.  A similar result holds for functors.

 \begin{prop}[Fisher-Palmquist, Newell, \cite{fisheradjoint}]For any functor $F\:\A\To \ab$, there is a natural transformation $$F\otimes \blank\To (F^*,\blank)$$ \end{prop}

\begin{cor}For any functor $F$, there is a natural transformation $$F\To F^{**}$$ obtained by evaluating the map $F\otimes \blank\To (F^*,\blank)$ at the bifunctor $\hom$.\end{cor}

From the defect sequence \begin{center}\begin{tikzpicture}
\matrix(m)[ampersand replacement=\&, matrix of math nodes, row sep=3em, column sep=2.5em, text height=1.5ex, text depth=0.25ex]
{0\&F_0\&F\&\big(w(F),\blank\big)\&F_1\&0\\}; 
\path[->, font=\scriptsize]
(m-1-1) edge node[auto]{}(m-1-2)
(m-1-2) edge node[auto]{}(m-1-3)
(m-1-3) edge node[auto]{$\varphi$}(m-1-4)
(m-1-4) edge node[auto]{}(m-1-5)
(m-1-5) edge node[auto]{}(m-1-6);
\end{tikzpicture}\end{center} and the fact that for any $F\in \fp(\A,\ab)$ we have $F^{**}\cong R^0F\cong \big(w(F),\blank\big)$, we have the following exact sequence

\begin{center}\begin{tikzpicture}
\matrix(m)[ampersand replacement=\&, matrix of math nodes, row sep=3em, column sep=2.5em, text height=1.5ex, text depth=0.25ex]
{0\&F_0\&F\&F^{**}\&F_1\&0\\}; 
\path[->, font=\scriptsize]
(m-1-1) edge node[auto]{}(m-1-2)
(m-1-2) edge node[auto]{}(m-1-3)
(m-1-3) edge node[auto]{$\varphi$}(m-1-4)
(m-1-4) edge node[auto]{}(m-1-5)
(m-1-5) edge node[auto]{}(m-1-6);
\end{tikzpicture}\end{center}  It is easily seen that the map $\varphi\:F\To F^{**}$ is precisely the map obtained by evaluating the natural transformation $$F\otimes\blank\To (F^*,\blank)$$ at the bifunctor $\hom$.  This allows one to calculate $F_0$ and $F_1$ explicitly.  Namely $F_k\cong \ext^{k+1}(\tr(F),\hom)$  where $$\ext^{k+1}\big(\tr(F),\hom\big)(X):=\ext^{k+1}\big(\tr(F),(\blank,X)\big)$$  

The proof is essentially the same as that given for Proposition 6.3 in \cite{coh}.    The argument given there is in terms of modules and it translates easily to $\fp(\A,\ab)$.  We can now completely explain the connection between the injective stabilization sequence of a finitely presented functor, the defect sequence, and the double dual sequence.  They are in fact all isomorphic.  That is, for any finitely presented functor $F$, there is a commutative diagram with exact rows

\begin{center}\begin{tikzpicture}
\matrix(m)[ampersand replacement=\&, matrix of math nodes, row sep=3em, column sep=1.5em, text height=1.5ex, text depth=0.25ex]
{0\&\ext^1(\tr(F),\hom)\&F\&F^{**}\&\ext^2(\tr(F),\hom)\&0\\
0\&F_0\&F\&\big(w(F),\blank\big)\&F_1\&0\\
0\&\overline{F}\&F\&R^0F\&\tilde{F}\&0\\}; 
\path[->, font=\scriptsize]
(m-1-2) edge node[auto]{$\cong$}(m-2-2)
(m-2-2) edge node[auto]{$\cong$}(m-3-2)
(m-1-3) edge node[auto]{$1_F$}(m-2-3)
(m-2-3) edge node[auto]{$1_F$}(m-3-3)
(m-1-4) edge node[auto]{$\cong$}(m-2-4)
(m-2-4) edge node[auto]{$\cong$}(m-3-4)
(m-1-5) edge node[auto]{$\cong$}(m-2-5)
(m-2-5) edge node[auto]{$\cong$}(m-3-5);
\path[->, font=\scriptsize]
(m-1-1) edge node[auto]{}(m-1-2)
(m-1-2) edge node[auto]{}(m-1-3)
(m-1-3) edge node[auto]{}(m-1-4)
(m-1-4) edge node[auto]{}(m-1-5)
(m-1-5) edge node[auto]{}(m-1-6)
(m-2-1) edge node[auto]{}(m-2-2)
(m-2-2) edge node[auto]{}(m-2-3)
(m-2-3) edge node[auto]{$\varphi$}(m-2-4)
(m-2-4) edge node[auto]{}(m-2-5)
(m-2-5) edge node[auto]{}(m-2-6)
(m-3-1) edge node[auto]{}(m-3-2)
(m-3-2) edge node[auto]{}(m-3-3)
(m-3-3) edge node[auto]{}(m-3-4)
(m-3-4) edge node[auto]{}(m-3-5)
(m-3-5) edge node[auto]{}(m-3-6);
\end{tikzpicture}\end{center}

From this discussion we can now state the \fp-injective stabilization formula which applies to any finitely presented functor.
\begin{thm}[\fp-Injective Stabilization Formula]For any $F\in \fp(\A,\ab)$ 
\begin{enumerate}
\item $F_0\cong\overline{F}\cong\ext^1(\tr(F),\hom)$
\item $F_1\cong \tr\u\tr(F)\cong\ext^2(\tr(F),\hom)$ 
\item $\u\tr\u\tr F$ is the image of the unit of adjunction $\varphi=u_F\:F\to F^{**}$
\end{enumerate} \end{thm}

\begin{cor}Suppose that $$0\To F\To G\To H\To 0$$ is a short exact sequence of finitely presented functors.  Then there is an exact sequence $$0\To F_0\To G_0\To H_0\To F_1\To G_1\To H_1$$\end{cor}

\begin{pf}The functor $R^0$ is the composition of the exact functor $w$ followed by the left exact Yoneda embedding.  Therefore it is a left exact functor.   As a result, the short exact sequence $$0\To F\To G\To H\To 0$$ embeds into the following commutative diagram with exact rows and columns:

\begin{center}\begin{tikzpicture}
\matrix(m)[ampersand replacement=\&, matrix of math nodes, row sep=3em, column sep=2.5em, text height=1.5ex, text depth=0.25ex]
{\&0\&0\&0\&\\
0\&F_0\&G_0\&H_0\&\\
0\&F\&G\&H\&0\\
0\&R^0F\&R^0G\&R^0H\&\\
\&F_1\&G_1\&H_1\&\\
\&0\&0\&0\\}; 
\path[->, font=\scriptsize]
(m-3-1) edge node[auto]{}(m-3-2)
(m-3-2) edge node[auto]{}(m-3-3)
(m-3-3) edge node[auto]{}(m-3-4)
(m-3-4) edge node[auto]{}(m-3-5)
(m-2-1) edge node[auto]{}(m-2-2)
(m-2-2) edge node[auto]{}(m-2-3)
(m-2-3) edge node[left]{}(m-2-4)
(m-4-1) edge node[below]{}(m-4-2)
(m-4-2) edge node[below]{}(m-4-3)
\testarr{4}{3}{4}{4}
\testarr 5 2 5 3
\testarr 5 3 5 4
\testarr  1 2 2 2
\testarr 1 3 2 3
\testarr 1 4 2 4
\testarr  2 2 3 2
\testarr 2 3 3 3
\testarr 2 4 3 4
\testarr  3 2 4 2
\testarr 3 3 4 3
\testarr 3 4 4 4
\testarr  4 2 5 2
\testarr 4 3 5 3
\testarr 4 4 5 4
\testarr  5 2 6 2
\testarr 5 3 6 3
\testarr 5 4 6 4;
\end{tikzpicture}\end{center}

Applying the snake lemma yields the exact sequence $$0\To F_0\To G_0\To H_0\To F_1\To G_1\To H_1\qquad\qed$$
\end{pf}

\section{\label{section6}Derived Functors of the Dual}

The category of finitely presented functors $\fp(\A^{op},\ab)$ has enough injectives under the assumption that $\A$ is a general abelian category with enough injectives.  This was first shown by Gentle in \cite{rongentle}.  The following is the approach employed there.   Define the category $\A/\mathcal{L}$ as follows.  The objects are left exact sequences $0\to A\to B\to C$ in $\A$.  The morphisms between two left exact sequences are the normal chain maps modulo those which factor through split left exact sequences.   After defining the category, Gentle comments that $\A/\mathcal{L}$ is equivalent to $\fp(\A^{op},\ab)$.  

\begin{prop}[Gentle, \cite{rongentle}, Prop 1.4]Left exact sequences of the form $0\to K\to I_1\to I_2$ with $I_1,I_2$ injectives in $\A$ are injective objects in $\A/\mathcal{L}$.\end{prop}  

Immediately after proving this proposition, Gentle shows that $\A/\mathcal{L}$ has enough injectives and injective dimension less than or equal to 2.   Since $\fp(\A^{op},\ab)$ is equivalent to $\A/\mathcal{L}$, this establishes that $\fp(\A^{op},\ab)$ has enough injectives whenever $\A$ has enough injectives.  This implies that if $\A$ has enough projectives, then $\fp(\A,\ab)$ which is equivalent to $\fp((\A^{op}) ^{op},\ab)$ has enough injectives.

\begin{prop}[Gentle, \cite{rongentle}]\label{enoughinj}If $\A$ has enough projectives, then $\fp(\A,\ab)$ has enough injectives and for every $F\in\fp(\A,\ab)$, there exists injective resolution
$$0\to F\to I^0\to I^1\to I^2\to 0$$ \end{prop}

\begin{ex}Let $R$ be any ring.  The category $\Mod(R)$ has enough injectives and enough projectives.  As an immediate application of Gentle's result, the categories $\fp(\Mod(R),\ab)$ and $\fp(\Mod(R)^{op},\ab)$ have enough injectives. \end{ex}

\begin{prop}[Auslander, \cite{coh}, Lemma 5.1]A functor $F\in \fp(\A,\ab)$ is injective if and only if $F$ is right exact.\end{prop}  

From this characterization, Gentle's result, and simple homological methods, one can establish that if $\A$ has enough projectives, then the injectives of $\fp(\A,\ab)$ are precisely the functors $F$ for which there is an exact sequence $$(Q,\blank)\To (P,\blank)\To F\To 0$$ with both $P$ and $Q$ projective.  We now return to the \fp-dual formula.  The functors $(\blank)^*$ defined by Fisher-Palmquist and Newell are left exact contravariant functors:

\begin{center}\begin{tikzpicture}
\matrix (m) [ampersand replacement= \&,matrix of math nodes, row sep=3em, 
column sep=6em,text height=1.5ex,text depth=0.25ex] 
{\fp(\Mod(R),\ab)\&\fp(\Mod(R)^{op},\ab)\\}; 
\path[->,thick, font=\scriptsize]
(m-1-1) edge[bend left] node[auto]{$(\blank)^*$}(m-1-2)
(m-1-2) edge[bend left] node[below]{$(\blank)^*$}(m-1-1);
\end{tikzpicture}\end{center}   Because the categories involved have enough injectives, the left derived functors $L^k(\blank)^*$ can be computed by using injective resolutions.

Let $F$ be a finitely presented functor and $$0\To F\To I^0\To I^1\To I^2\To 0$$ its injective resolution.   The functor $L^k(\blank)^*$ evaluated at $F$ is the $k$-th homology of the complex $$0\To (I^2)^*\To (I^1)^*\To (I^0)^*\To 0$$ It is easily verified that for all $k\ge 1$, $L^k(\blank)^*\cong 0$ because $(\blank)^*$ is left exact.     Set $n=L^0(\blank)^*$.  Then the functor $n(F)$ is determined by the exact sequence $$(I^1)^*\To (I^0)^*\To n(F)\To0$$   In fact, because the source categories have injective dimension less than or equal to 2 and because $(\blank)^*$ is left exact, it follows that the functors  

\begin{center}\begin{tikzpicture}
\matrix (m) [ampersand replacement= \&,matrix of math nodes, row sep=3em, 
column sep=6em,text height=1.5ex,text depth=0.25ex] 
{\fp(\Mod(R),\ab)\&\fp(\Mod(R)^{op},\ab)\\}; 
\path[->,thick, font=\scriptsize]
(m-1-1) edge[bend left] node[auto]{$n=L^0(\blank)^*$}(m-1-2)
(m-1-2) edge[bend left] node[below]{$n=L^0(\blank)^*$}(m-1-1);
\end{tikzpicture}\end{center} are exact.  
Since $n$ agrees with $(\blank)^*$ on injectives, for any finitely presented injective functor $I$, we have $n(I)\cong I^*\cong (\blank, w(I))$.
With this notation, we now show the following.

\begin{prop}If $n(F)$ is representable, then $n(F)\cong (\blank,w(F))\cong F^*$.  In particular, if $n(F)=0$, then $w(F)=0$.  \end{prop}

\begin{pf}Suppose that $n(F)\cong(\blank, A)$.  Take any injective resolution $$0\To F\To I^0\To I^1\To I^2\To 0$$ and apply the exact functor $n$ yielding commutative diagrams with exact rows \dia{0\ar[r]&\big(\blank, w(I^2)\big)\ar@{=}[d]\ar[r]&\big(\blank, w(I^1)\big)\ar[r]\ar@{=}[d]&\big(\blank, w(I^0)\big)\ar[r]\ar@{=}[d]& n(F)\ar[r]\ar[d]^{\cong}&0\\
0\ar[r]&\big(\blank, w(I^2)\big)\ar[r]&\big(\blank, w(I^1)\big)\ar[r]&\big(\blank, w(I^0)\big)\ar[r]&(\blank, A)\ar[r]&0} Evaluating the bottom row at $R$ and using the natural isomorphism $(R,X)\cong X$ yields exact sequence $$0\To w(I^2)\To w(I^1)\To w(I^0)\To A\To 0$$  Since $w$ is exact, it follows that by applying $w$ to the original injective resolution of $F$, we have a commutative diagram with exact rows \dia{0\ar[r]& w(I^2)\ar@{=}[d]\ar[r]&w(I^1)\ar[r]\ar@{=}[d]& w(I^0)\ar[r]\ar@{=}[d]& w(F)\ar[r]\ar[d]^{\cong}&0\\
0\ar[r]&w(I^2)\ar[r]& w(I^1)\ar[r]&w(I^0)\ar[r]&A\ar[r]&0}  Therefore $A\cong w(F)$ and hence $n(F)\cong \big(\blank, w(F)\big)$.  Now if $n(F)=0$, then $n(F)\cong(\blank,w(F))=0$ and by Yoneda's lemma, $w(F)=0$.  $\qed$
\end{pf}  

This establishes that if $n(F)=0$, then $F$ arises from a short exact sequence in $\Mod(R)$ since $w(F)=0$.  In particular, because $n$ is exact, its kernel is a Serre subcategory of $\fp(\Mod(R),\ab)$ and actually a Serre subcategory of $\fp_0(\Mod(R),\ab)$.

\begin{ex}Consider any finite non-trivial abelian group $X$.  Because $X$ is a torsion module, it cannot map into $\ZZ$ and so $X^*\cong 0$.  By Corollary \ref{tdf},   $(X\otimes \blank)^*\cong (\blank, X^*)\cong 0$.  Since $X\otimes \blank$ is right exact, it  is an injective object in $\fp(\Mod(\ZZ^{op}),\ab)$. Therefore $n(X\otimes \blank)\cong\big(X\otimes\blank\big)^*\cong0$ while $X\otimes \blank$ is not a trivial functor as $(X\otimes \ZZ)\cong X\ne 0$.  Hence $n$ vanishes on non-trivial functors. \end{ex}

Because every functor $F\in \fp(\Mod(R),\ab)$ has an injective resolution $$0\To F\To I^0\To I^1\To I^2\To 0$$ and $n$ agrees with $(\blank)^*$ on injectives we have a commutative diagram with exact rows \dia{0\ar[r] &(I^2)^*\ar[r]\ar[d]_{\cong}&(I^1)^*\ar[r]\ar[d]_{\cong}&(I^0)^*\ar[r]\ar[d]_{\cong}&n(F)\ar[r]\ar[d]_{\cong}&0\\
0\ar[r]&(\blank, w(I^2))\ar[r]&(\blank, w(I^1))\ar[r]&(\blank, w(I^0))\ar[r]&n(F)\ar[r]&0}  If $n(F)=0$, then since $(\blank, w(I^0))$ is projective in $\fp(\Mod(R)^{op},\ab)$ the sequence $$0\To (\blank, w(I^2))\To (\blank, w(I^1))\To (\blank, w(I^0))\To 0$$ splits.  Evaluating at $R$ gives the split exact sequence $$0\To w(I^2)\To w(I^1)\To w(I^0)\To 0$$ On the other hand, if the sequence $$0\To w(I^2)\To w(I^1)\To w(I^0)\To 0$$ is split, then the sequence $$0\To (\blank, w(I^2))\To (\blank, w(I^1))\To (\blank, w(I^0))\To 0$$ is also split and hence $n(F)=0$.   This establishes that $n(F)=0$ if and only if given any injective resolution $0\to F\to I^0\to I^1\to I^2\to 0$, the sequence $0\to w(I^2)\to w(I^1)\to w(I^0)\to 0$ is a split exact sequence.  

\begin{ex}Let $X$ be any module and take a projective resolution $$\cdots\To P_3\To P_2\To P_1\To P_0\To X\To 0$$  We have a commutative diagram with exact rows and columns

\dia{&&&0\ar[d]&0\ar[d]&\\
0\ar[r]&(X,\blank)\ar[r]\ar@{=}[d]&(P_0,\blank)\ar[r]\ar@{=}[d]&(\u X,\blank)\ar[r]\ar[d]&\ext^1(X,\blank)\ar[r]\ar[d]&0\\
0\ar[r]&(X,\blank)\ar[r]&(P_0,\blank)\ar[r]&(P_1,\blank)\ar[r]\ar[d]&I^0\ar[r]\ar@{-->}[dl]&0\\
&&&(P_2,\blank)
}  and a simple diagram chase shows that $$0\To \ext^1(X,\blank)\To I^0\To (P_2,\blank)$$ is exact.  Moreover, $I^0$ and $I^1=(P_2,\blank)$ are injectives in $\fp(\Mod(R),\ab)$ because they are both right exact and hence  $0\to \ext^1(X,\blank)\to I^0\to I^1$ is an injective copresentation of $\ext^1(X,\blank)$.    Since $\ext^1(X,\blank)$ arises from a short exact sequence, it follows that $w(\ext^1(X,\blank))=0$.   However, from the exact sequence $$0\To w(I^0)\To P_1\To P_0\To X\To 0$$ we have $w(I^0)\cong \u^2 X$.   From the commutative diagram with exact rows 
\dia{w(I^1)\ar[d]_{\cong}\ar[r]& w(I^0)\ar[r]\ar[d]_{\cong} &w(\ext^1(X,\blank))\ar[d]_{\cong}\ar[r]& 0\\
P_2\ar[r]&\u^2 X\ar[r]&0} it follows that $w(I^1)\to w(I^0)\to 0$ splits if and only if   $P_2\to \u^2X\to 0$ splits, which occurs if and only if $\u^2X$ is a projective module.  Therefore, by the discussion immediately preceding this example, if $X$ has projective dimension larger than 2, then  $n\big(\ext^1(X,\blank)\big)\ne 0$ and yet $w\big(\ext^1(X,\blank)\big)=0$.   \end{ex}

The means that the functors $n$ and $w$ in general will have different kernels  resulting in an open question.  What are the short exact sequences in $\Mod(R)$ that correspond to the finitely presented functors $F$ for which $n(F)=0$? It seems feasible that these sequences are worth investigating though at this time we do not have a clear description.   We conclude by combining the results from above into the following theorem.

\begin{thm}The exact functor $n=L^0(\blank)^*$ is the only non-trivial left derived functor of the dual. The kernel of $n$ is a Serre subcategory of $\fp(\Mod(R),\ab)$ which vanishes on functors whose presentations arise from certain short exact sequences.  For any $F\in \fp(\Mod(R),\ab)$ the following are equivalent:
\begin{enumerate}
\item $n(F)=0$.
\item For every injective resolution $$0\To F\To I^0\To I^1\To I^2\To 0$$ the sequence $$0\To w(I^2)\To w(I^1)\To w(I^0)\To 0$$ is a split exact sequence.
\item $w(F)=0$ and for any injective copresentation $$0\To F\To I^0\To I^1$$ the map $$w(I^1)\To w(I^0)$$ is a retraction.
\end{enumerate}
\end{thm}

\section{\label{section7}Detecting Almost Split Sequences}
 
The purpose of this section is to demonstrate how studying the category of finitely presented functors $\fp(\A,\ab)$ and using functorial techniques can reveal information about the category $\A$ which, one may argue, is more easily seen in the category $\fp(\A,\ab)$.  For example, the defect $w$ predicts existence of special short exact sequences in $\A$ under the right conditions.  These are the so called almost split sequences discovered and studied in depth by Auslander and Reiten.  The proof of their existence first appears in \cite{ar3}.   In this section, we show how to recover almost split sequences using the defect.

\begin{defn}A functor $S\:\A\to \ab$ is called simple if $S\ne 0$ and any non-zero morphism $F\to S$ is an epimorphism.\end{defn}

\begin{prop}\label{simpledefect}Suppose that $S\:\A\To \ab$ is a functor satisfying the following:
\begin{enumerate}
\item $S$ is a finitely presented simple functor.
\item $w(S)=0$.
\item $S$ has a minimal projective presentation.
\end{enumerate}
Then there exists a short exact sequence $$0\To X\To Y\To Z\To 0$$ such that given any $u$ \begin{center}\begin{tikzpicture}
\matrix(m)[ampersand replacement=\&, matrix of math nodes, row sep=3em, column sep=2.5em, text height=1.5ex, text depth=0.25ex]
{0\&X\&Y\&Z\&0\\
\&K\&\&\\}; 
\path[->]
(m-1-1) edge node[auto]{}(m-1-2)
(m-1-2) edge node[auto]{$f$}(m-1-3)
(m-1-3) edge node[auto]{}(m-1-4)
(m-1-4) edge node[auto]{}(m-1-5)
(m-1-2) edge node[left]{$u$}(m-2-2);
\end{tikzpicture}\end{center}
either $u$ is a section or there exists $g\:Y\To K$ such that the following diagram commutes

\begin{center}\begin{tikzpicture}
\matrix(m)[ampersand replacement=\&, matrix of math nodes, row sep=3em, column sep=2.5em, text height=1.5ex, text depth=0.25ex]
{0\&X\&Y\&Z\&0\\
\&K\&\&\\}; 
\path[->]
(m-1-1) edge node[auto]{}(m-1-2)
(m-1-2) edge node[auto]{$f$}(m-1-3)
(m-1-3) edge node[auto]{}(m-1-4)
(m-1-4) edge node[auto]{}(m-1-5)
(m-1-2) edge node[left]{$u$}(m-2-2)
(m-1-3)edge node[below]{$g$}(m-2-2);
\end{tikzpicture}\end{center}

\end{prop}

\begin{pf}Take a minimal presentation of $S$ $$0\To (Z,\blank)\To (Y,\blank)\To (X,\blank)\To S\To 0$$  Since $w(S)=0$, the following sequence is a short exact exact sequence 
$$0\To X\To Y\To Z\To 0$$  
Take any morphism $u\:X\To K$.  
\begin{center}\begin{tikzpicture}
\matrix(m)[ampersand replacement=\&, matrix of math nodes, row sep=3em, column sep=2.5em, text height=1.5ex, text depth=0.25ex]
{0\&X\&Y\&Z\&0\\
\&K\&\&\\}; 
\path[->]
(m-1-1) edge node[auto]{}(m-1-2)
(m-1-2) edge node[auto]{$f$}(m-1-3)
(m-1-3) edge node[auto]{}(m-1-4)
(m-1-4) edge node[auto]{}(m-1-5)
(m-1-2) edge node[left]{$u$}(m-2-2);
\end{tikzpicture}\end{center}

Taking the pushout of this diagram results in the following commutative diagram with exact rows.

\begin{center}\begin{tikzpicture}
\matrix(m)[ampersand replacement=\&, matrix of math nodes, row sep=3em, column sep=2.5em, text height=1.5ex, text depth=0.25ex]
{0\&X\&Y\&Z\&0\\
0\&K\&E\&Z\&0\\
\&\&\&\\
\&\&\\}; 
\path[->]
(m-1-1) edge node[auto]{}(m-1-2)
(m-1-2) edge node[auto]{$f$}(m-1-3)
(m-1-3) edge node[auto]{}(m-1-4)
(m-1-4) edge node[auto]{}(m-1-5)
(m-1-2) edge node[left]{$u$}(m-2-2)
(m-2-1) edge node[auto]{}(m-2-2)
(m-2-2) edge node[auto]{}(m-2-3)
(m-2-3) edge node[auto]{}(m-2-4)
(m-2-4) edge node[auto]{}(m-2-5)
(m-1-3) edge node[auto]{}(m-2-3)
(m-1-4) edge node[auto]{$1$}(m-2-4);
\end{tikzpicture}\end{center}

Using the Yoneda embedding and extending to include cokernels where necessary, we obtain a commutative diagram with exact rows:

\begin{center}\begin{tikzpicture}
\matrix(m)[ampersand replacement=\&, matrix of math nodes, row sep=3em, column sep=2.5em, text height=1.5ex, text depth=0.25ex]
{0\&(Z,\blank)\&(E,\blank)\&(K,\blank)\&F\&0\\
0\&(Z,\blank)\&(Y,\blank)\&(X,\blank)\&S\&0\\}; 
\path[->]
(m-1-1) edge node[auto]{}(m-1-2)
(m-1-2) edge node[auto]{}(m-1-3)
(m-1-3) edge node[auto]{}(m-1-4)
(m-1-4) edge node[auto]{}(m-1-5)
(m-1-2) edge node[left]{$1$}(m-2-2)
(m-2-1) edge node[auto]{}(m-2-2)
(m-2-2) edge node[auto]{}(m-2-3)
(m-2-3) edge node[below]{$(f,\blank)$}(m-2-4)
(m-2-4) edge node[auto]{}(m-2-5)
(m-1-3) edge node[auto]{}(m-2-3)
(m-1-4) edge node[auto]{$(u,\blank)$}(m-2-4)
(m-1-5) edge node[auto]{$\alpha$}(m-2-5)
(m-2-5)edge node[auto]{}(m-2-6)
(m-1-5)edge node[auto]{}(m-1-6);
\end{tikzpicture}\end{center}

By assumption, $S$ is simple and therefore, either $\alpha$ is an epimorphism or $\alpha=0$.

\dbf{Case 1: $\alpha=0$}
 
 If $\alpha=0$, there exists $(g,\blank)\:(K,\blank)\To (Y,\blank)$ such that 
 $$(f,\blank)(g,\blank)=(u,\blank)$$ 

\begin{center}\begin{tikzpicture}
\matrix(m)[ampersand replacement=\&, matrix of math nodes, row sep=3em, column sep=2.5em, text height=1.5ex, text depth=0.25ex]
{0\&(Z,\blank)\&(E,\blank)\&(K,\blank)\&F\&0\\
0\&(Z,\blank)\&(Y,\blank)\&(X,\blank)\&S\&0\\}; 
\path[->]
(m-1-1) edge node[auto]{}(m-1-2)
(m-1-2) edge node[auto]{}(m-1-3)
(m-1-3) edge node[auto]{}(m-1-4)
(m-1-4) edge node[auto]{}(m-1-5)
(m-1-2) edge node[left]{$1$}(m-2-2)
(m-2-1) edge node[auto]{}(m-2-2)
(m-2-2) edge node[auto]{}(m-2-3)
(m-2-3) edge node[below]{$(f,\blank)$}(m-2-4)
(m-2-4) edge node[auto]{}(m-2-5)
(m-1-3) edge node[auto]{}(m-2-3)
(m-1-4) edge node[auto]{$(u,\blank)$}(m-2-4)
(m-1-5) edge node[auto]{$\alpha=0$}(m-2-5)
(m-2-5)edge node[auto]{}(m-2-6)
(m-1-5)edge node[auto]{}(m-1-6)
(m-1-4)edge node[above]{}(m-2-3);
\end{tikzpicture}\end{center}

which results in the commutative diagram

\begin{center}\begin{tikzpicture}
\matrix(m)[ampersand replacement=\&, matrix of math nodes, row sep=3em, column sep=2.5em, text height=1.5ex, text depth=0.25ex]
{0\&X\&Y\&Z\&0\\
\&K\&\&\\}; 
\path[->]
(m-1-1) edge node[auto]{}(m-1-2)
(m-1-2) edge node[auto]{$f$}(m-1-3)
(m-1-3) edge node[auto]{}(m-1-4)
(m-1-4) edge node[auto]{}(m-1-5)
(m-1-2) edge node[left]{$u$}(m-2-2)
(m-1-3)edge node[below]{$g$}(m-2-2);
\end{tikzpicture}\end{center}

\dbf{Case 2: $\alpha$ is an epimorphism}

In this case the composition is $(K,\blank)\To F\To S$ an epimorphism but then the composition $$(K,\blank)\overset{(u,\blank)}{\To} (X,\blank)\To S$$ is an epimorphism.  But $(X,\blank)\To S$ is a projective cover and therefore, $(u,\blank)$ must be an epimorphism.    Hence $u$ is a section.  $\qed$
\end{pf}

We conclude with an example justifying studying the defect more closely.  The fact that $w$ vanishes on those functors whose presentations arise from short exact sequences combined with the fact that almost split sequences are specific examples of short exact sequences indicated that there may be a connection between the two.  Since for an abelian category $\A$, the Yoneda embedding $$\Y\:\A\To \fp(\A,\ab)$$ is not exact, the finitely presented functors can be thought of as a type of ``noise'' created by applying the non-exact embedding.   Keeping with this analogy, the ``noise" created by the almost spit sequences in the functor category, one can argue, makes them easier to ``hear''.  While the correspondence between the simple functors and the almost split sequences is well known, it is interesting to find these sequences using the defect.

\begin{ex}Let $\Lambda$ be a finite dimensional $k$-algebra.  Then $\mod(\Lambda)$ is an abelian category with enough projectives and injectives.  In \cite{fart}, Auslander shows that $\fp(\mod(\Lambda),\ab)$ has minimal projective presentations.  Moreover, Auslander classified the simple functors as follows.  $S\:\mod(\Lambda)\To \ab$ is simple if and only if there exists a unique indecomposable $\Lambda$-module  $N$ such that $S(N)\ne 0$.  This establishes a bijection between indecomposable objects in $\mod(\Lambda)$ and simple functors.  We will say that the simple functor $S$ is determined by $N$ if $N$ is the unique indecomposable module where $S$ does not vanish.   Auslander further shows in \cite{fart} that all simple functors $S$ are finitely presented and hence have minimal projective resolutions.  It is easily seen that $w(S)=0$ when $S$ is determined by a non-injective $\Lambda$-module $N$.  Therefore, by Proposition \ref{simpledefect} every simple functor determined by non-injective $N$ gives rise to an almost split sequence.  \end{ex}

\bibliographystyle{amsplain}
\bibliography{slayer}

\end{document}